\documentclass[reqno]{amsart}
\usepackage{pictex}
\usepackage{mathrsfs}
\usepackage[colorlinks]{hyperref}
\usepackage{color}
\begin{document}

%
%

\def\labelenumi{(\theenumi)}

\newtheorem{thm}{Theorem}[section]
\newtheorem{lem}[thm]{Lemma}
\newtheorem{cor}[thm]{Corollary}
\newtheorem{add}[thm]{Addendum}
\newtheorem{prop}[thm]{Proposition}
\theoremstyle{definition}
\newtheorem{defn}[thm]{Definition}
\theoremstyle{remark}
\newtheorem{rmk}[thm]{Remark}
\newtheorem{example}[thm]{{\bf Example}}

\newcommand{\cR}{\,\,\not \!\! R}
\newcommand{\SLtwoC}{\mathrm{SL}(2,{\mathbf C})}
\newcommand{\SLtwoR}{\mathrm{SL}(2,{\mathbb R})}
\newcommand{\PSLtwoC}{\mathrm{PSL}(2,{\mathbf C})}
\newcommand{\PSLtwoR}{\mathrm{PSL}(2,{\mathbb R})}
\newcommand{\SLtwoZ}{\mathrm{SL}(2,{\mathbb Z})}
\newcommand{\PSLtwoZ}{\mathrm{PSL}(2,{\mathbb Z})}
\newcommand{\CmodTwoPiIZ}{{\mathbf C}/2\pi i {\mathbb Z}}

\newcommand{\Cnozero}{{\mathbf C}\backslash \{0\}}
\newcommand{\Cinfty}{{\mathbf C}_{\infty}}
\newcommand{\HH}{{\mathbb H}^2}
\newcommand{\HHH}{{\mathbb H}^3}

\newcommand{\OmegaH}{\Omega/\langle H \rangle}
\newcommand{\hatOmegaHstar}{\hat \Omega/\langle H_{\ast}\rangle}
\newcommand{\SurfG}{\Sigma_g}
\newcommand{\TriangG}{T_g}
\newcommand{\TriangGOne}{T_{g,1}}
\newcommand{\ProjG}{\mathcal{P}_g}
\newcommand{\TeichG}{\mathcal{T}_g}
\newcommand{\CirclePackGTau}{\mathsf{CPS}_{g,\tau}}
\newcommand{\CrossRatio}{{\bf c}}
\newcommand{\CrossRatioGTau}{\mathcal{C}_{g,\tau}}
\newcommand{\CrossRatioOneTau}{\mathcal{C}_{1,\tau}}
\newcommand{\DeformGTau}{\mathcal{C}_{g,\tau}}
\newcommand{\Forget}{\mathit{forg}}
\newcommand{\Uniform}{\mathit{u}}
\newcommand{\MCSS}{{\mathcal S}^{\rm mc}_{\rm alg}}
\newcommand{\Section}{\mathit{sect}}
\newcommand{\SLTwoC}{\mathrm{SL}(2,\mathbf{C})}
\newcommand{\SLTwoR}{\mathrm{SL}(2,\mathbf{R})}
\newcommand{\PSLTwoC}{\mathrm{PSL}(2,\mathbf{C})}
\newcommand{\GLTwoC}{\mathrm{GL}(2,\mathbf{C})}
\newcommand{\PSLTwoR}{\mathrm{PSL}(2,\mathbf{R})}
\newcommand{\PSLTwoZ}{\mathrm{PSL}(2,\mathbf{Z})}
\newcommand{\SLTwoZ}{\mathrm{SL}(2,\mathbf{Z})}
\newcommand{\nnn}{\noindent}
\newcommand{\MCG}{{\mathcal {MCG}}}
\newcommand{\MMap}{{\bf \Phi}_{\mu}}
\def\square{\hfill${\vcenter{\vbox{\hrule height.4pt \hbox{\vrule width.4pt
height7pt \kern7pt \vrule width.4pt} \hrule height.4pt}}}$}

\newenvironment{pf}{\noindent {\sl Proof.}\quad}{\square \vskip 12pt}

\title{McShane's identity
for classical Schottky groups}
\author{Ser Peow Tan, Yan Loi Wong, and Ying Zhang}
\address{Department of Mathematics \\ National University of Singapore \\
2 Science Drive 2 \\ Singapore 117543} \email{mattansp@nus.edu.sg;
matwyl@nus.edu.sg; scip1101@nus.edu.sg}
\address{and the third author}
\address{Department of Mathematics \\Yangzhou University \\Yangzhou 225002 \\P. R. China}
\email{yingzhang@yzu.edu.cn}

\thanks{The authors are partially supported by the National University
of Singapore academic research grant R-146-000-056-112. The third
author is also partially supported by the National Key Basic
Research Fund (China) G1999075104.}

%
%

\begin{abstract}
In \cite{mcshane1998im}, Greg McShane demonstrated a  remarkable
identity for the lengths of simple closed geodesics on cusped
hyperbolic surfaces. This was generalized by the authors in
\cite{tan-wong-zhang2004cone-surfaces} to hyperbolic
cone-surfaces, possibly with cusps and/or geodesic boundary. In
this paper, we generalize the identity further to the case of
classical Schottky groups. As a consequence, we obtain some
surprising new identities in the case of fuchsian Schottky groups.
For classical Schottky groups of rank 2, we also give
generalizations of the Weierstrass identities, given by McShane in
\cite{mcshane2004blms}.
\end{abstract}

\maketitle

\section{{ Introduction}}\label{s:intro}
\vskip 5pt

In \cite{mcshane1998im} Greg McShane demonstrated a striking
identity for the lengths of simple closed geodesics on cusped
hyperbolic surfaces. In \cite{tan-wong-zhang2004cone-surfaces}, we
gave a generalization of McShane's identity to hyperbolic
cone-surfaces possibly with cusps and/or geodesic boundary (a
version for compact hyperbolic surfaces with geodesic boundary was
also given independently by M. Mirzakhani in
\cite{mirzakhani2004preprint}). To state the generalization, we
first define two functions $G(x,y,z)$ and $S(x,y,z)$ as follows.

\begin{defn}\label{def:GandS}
\begin{eqnarray}
G(x,y,z):=2
\tanh^{-1}\left(\frac{\sinh(x)}{\cosh(x)+\exp(y+z)}\right),
\end{eqnarray}
\begin{eqnarray}
S(x,y,z):=\tanh^{-1}\left(\frac{\sinh(x)\sinh(y)}{\cosh(z)+\cosh(x)\cosh(y)}\right).
\end{eqnarray}
\end{defn}

Note that here, for a complex number $x$, $\tanh^{-1}(x)$ is
defined to have imaginary part in $(-\pi/2, \pi/2]$, and that the
functions $G(x,y,z)$ and $S(x,y,z)$ are analytic in the arguments
$x,y$ and $z$.

\vskip 5pt

Recall, as defined in \cite{tan-wong-zhang2004cone-surfaces}, that
a {\it compact hyperbolic cone-surface} is a compact topological
surface, $M$, equipped with a hyperbolic cone structure so that
each boundary component of $M$ is a smooth geodesic, and there are
a finite number of points in the topological interior of $M$,
which form the set of  cusps and cone points. The {\it geometric
boundary} of $M$ consists of all the cusps, cone points and
boundary geodesics. The complement of the geometric boundary in
$M$ is its {\it geometric interior}. Recall also that (a) a
geometric boundary component of a compact hyperbolic cone-surface
$M$ is either a cusp, a cone point or a boundary geodesic; (b) an
interior generalized simple closed geodesic is either a simple
closed geodesic in the geometric interior of $M$ or a degenerate
one which is the double cover of a simple geodesic arc connecting
two angle $\pi$ cone points; (c) a generalized simple closed
geodesic is either an interior one as in (b) or a boundary one,
that is, a geometric boundary component.

The generalized identity in \cite{tan-wong-zhang2004cone-surfaces}
for compact hyperbolic cone-surfaces can be stated as follows.

\begin{thm}\label{thm:complexified non-cusp cases}{\rm(Theorem 9.3
\cite{tan-wong-zhang2004cone-surfaces})}\, For a compact
hyperbolic cone-surface $M$ with all cone angles in $(0, \pi]$,
let $\Delta_0, \Delta_1$, $\cdots, \Delta_m$ be its geometric
boundary components  with complex lengths $L_0, L_1, \cdots, L_m$
respectively. If $\Delta_0$ is a cone point or a boundary geodesic
then
\begin{eqnarray}\label{eqn:reform of non-cusp cases with GS}
\sum_{\alpha, \beta}
G\left(\frac{L_0}{2},\frac{|\alpha|}{2},\frac{|\beta|}{2}\right) +
\sum_{j=1}^{m}\sum_{\beta}
S\left(\frac{L_0}{2},\frac{L_j}{2},\frac{|\beta|}{2}\right)=
\frac{L_0}{2},
\end{eqnarray} where the first sum is taken over all
unordered pairs of generalized simple closed geodesics $\alpha,
\beta$ on $M$ such that $\alpha, \beta$ bound with $\Delta_0$ an
embedded pair of pants on $M$ {\rm(}note that one of $\alpha,
\beta$ might be a geometric boundary component\,{\rm)} and the
sub-sum in the second sum is taken over all interior simple closed
geodesics $\beta$ such that $\beta$ bounds with $\Delta_j$ and
$\Delta_0$ an embedded pair of pants on $M$. Furthermore, all the
infinite series in {\rm(\ref{eqn:reform of non-cusp cases with
GS})} converge absolutely.
\end{thm}

\vskip 10pt
\begin{rmk}\hfill
\begin{enumerate}

\item[(i)] $|\alpha|$ denotes the complex length of the generalized
simple closed geodesic $\alpha$, which is the usual hyperbolic
length if $\alpha$ is a geodesic, is equal to $i \theta$ if
$\alpha$ is a cone point of cone angle $\theta$ and is $0$ if
$\alpha$ is a cusp point. Similarly, $L_i$ is strictly positive,
pure imaginary or $0$ depending on whether $\Delta_i$ is a
geodesic boundary, cone point or cusp respectively.


\item[(ii)] As explained in
\cite{tan-wong-zhang2004cone-surfaces}, the above identity unifies
all the various known versions of McShane's identity in the real
case by using the complex lengths of the geometric boundary
components. A version for surfaces with geodesic boundary and no
cone points was also given by Mirzakhani in
\cite{mirzakhani2004preprint}. See also
\cite{akiyoshi-miyachi-sakuma2002preprint},
\cite{akiyoshi-miyachi-sakuma2004preprint},
\cite{bowditch1996blms}, \cite{bowditch1997t},
\cite{bowditch1998plms}, \cite{mcshane2004blms} and
\cite{sakuma1999sk} for other variations and generalizations of
the original identity.
\end{enumerate}
\end{rmk}

In this paper, motivated by the idea of considering the
complexification of the boundary lengths, we extend Theorem
\ref{thm:complexified non-cusp cases} to classical Schottky groups
by analytic continuation. This is possible because firstly, the
marked classical Schottky space, appropriately parametrized, is a
connected open subset of the parametrization space. Secondly, both
sides of the identity (\ref{eqn:reform of non-cusp cases with GS})
are analytic functions of the parameters for the marked classical
Schottky space. Finally, all the infinite series on the left hand
side of (\ref{eqn:reform of non-cusp cases with GS}) converge
absolutely in the marked classical Schottky space. However, we
will need to re-interpret $|\alpha|$ as the complex length of
$\alpha$ in the formula, and there are some subtleties involved as
we will need a specific choice of the half-lengths (recall that
the complex length is defined up to multiples of $2 \pi i$, and
that there are two possible choices for the half lengths, up to
multiples of $2 \pi i$). The exact statement (Theorem
\ref{thm:mcshane schottky}) requires a choice of a lift from a
representation of the free group $F_n$ into $\PSLTwoC$ to
$\SLTwoC$ (the exact choice of the lift is not important), a
reformulation of Theorem \ref{thm:complexified non-cusp cases} to
a more algebraic form, precise definitions for a marked classical
Schottky group, the marked classical Schottky space, as well as a
fuchsian marking for the marked classical Schottky space. The
basic idea is that a fuchsian marking corresponds to a hyperbolic
surface $M$ with geodesic boundary. Fixing a boundary component
$\Delta_0$ for $M$, we have that the identity (\ref{eqn:reform of
non-cusp cases with GS}) holds. Now as we deform away from the
fuchsian marking to an arbitrary point in the marked classical
Schottky space, the identity continues to hold by analytic
continuation, if we can show that the infinite series converge
absolutely and uniformly on compact subsets of the marked
classical Schottky space. In particular, if we deform to a
different fuchsian marked classical Schottky group, the original
identity still holds even though the corresponding hyperbolic
surface $M'$ may be of a different topological type from the
original surface $M$ (or it may be of the same topological type
but with a completely different marking). This gives new
identities for the surface $M'$. For example, in the case when the
Schottky group is of rank $2$, $M$ may be the hyperbolic one-holed
torus and $M'$ the hyperbolic three-holed sphere (pair of pants),
so that we obtain new identities for the hyperbolic pair of pants
with geodesic boundary, different from the trivial identity
obtained by a direct application of Theorem \ref{thm:complexified
non-cusp cases}.

The rest of this paper is organized as follows. In \S \ref{s:intro
schottky} we give some basic facts about classical Schottky groups
and give a precise definition of  marked classical Schottky space,
as well as a parametrization for the space. In \S \ref{s:mcshane
schottky} we state and prove the main result, Theorem
\ref{thm:mcshane schottky}. In \S \ref{s:Weierstrass} we state and
discuss the Weierstrass identity for marked rank 2 classical
Schottky groups (Theorem \ref{thm:Weierstrass}). Finally, in \S
\ref{s:example schottky}, we analyze an example to show that our
generalization of McShane's identity to classical Schottky groups
does give some surprising new identities for hyperbolic surfaces
with geodesic boundary.

\vskip 10pt

\noindent {\it Acknowledgements.} We would like to thank Caroline
Series, Bill Goldman,  Makoto Sakuma and Greg McShane for their
encouragement and/or helpful conversations.

\vskip 15pt
\section{Marked classical
Schottky groups}\label{s:intro schottky} \vskip 5pt

In this section we  state some basic facts about marked classical
Schottky groups. See \cite{marden1974}, \cite{maskit2000cm256},
\cite{maskit2002cm311}, \cite{button1998gtm1} and
\cite{button2003plms} for a more complete study of various
Schottky spaces. Note that the terminology is not completely
standardized, we use the terminology which is best suited for our
purposes, in particular, we need to pay special attention to the
marking in order to get a clean and precise statement of our
result. Hence, we define a marked classical Schottky group and
marked classical Schottky space in an analogous way to that of
marked hyperbolic structures and Teichm\"uller space.

\vskip 6pt Fix $n \ge 2$ and let $F_n=\langle a_1, \ldots,a_n
\rangle$ be a free group of rank $n$, where $\{a_1,\ldots, a_n\}$
is a (fixed) distinguished, ordered set of generators for $F_n$,
which will provide the marking. Let $\Cinfty$ be the extended
complex plane, which we also identify with the Riemann sphere, and
also the ideal boundary of $\HHH$.

\begin{defn}\label{def:markedschottkygroup}
A {\it {\rm(}marked\,{\rm)} classical Schottky group} (of rank
$n$) is a discrete, faithful representation $\rho:F_n \to
\PSLTwoC$ such that there is a region $D \subset \Cinfty $, where
$D$ is bounded by $2n$ disjoint geometric circles  $C_1, C'_1,
\cdots, C_n, C'_n$ in $\Cinfty$, so that, for $i=1, \ldots, n$,
$\rho(a_i)(C_i)= C'_i$, and $\rho(a_i)(D) \cap D = \emptyset$.
\end{defn}

Note that the circles $C_i,C_i'$ are not uniquely determined by
$\rho$. Also, $\rho(a_i)$ is strictly loxodromic, with an
attracting and repelling ideal fixed point, denoted by ${\rm
Fix}^{+}\rho(a_i)$ and ${\rm Fix}^{-}\rho(a_i)$ respectively. The
image $\Gamma:=\rho(F_n)$  is often referred to in the literature
as a classical Schottky group; it is a Schottky group if we only
require $C_i, C_i'$, $i=1, \ldots,n$ to be disjoint simple closed
curves. Two (marked) classical Schottky groups $\rho_1$ and
$\rho_2$ are equivalent if the representations are conjugate to
each other.

If $\rho$ is a marked classical Schottky group, we denote by
$\tilde \rho$ any lift of $\rho$ to  a representation into
$\SLTwoC$. Note that since $F_n$ is free, $\rho$ can always be
lifted and there are $2^n$ possible lifts.

\begin{rmk}{\rm It was shown by Marden \cite{marden1974} that there exist
non-classical Schottky groups for every $n \ge 2$. An explicit
example of a non-classical Schottky group was constructed by
Yamamoto \cite{yamamoto1991dmj}. On the other hand, Button
\cite{button1998gtm1} has proved that all fuchsian Schottky groups
are classical (but in general not on every set of
generators).}\end{rmk}

\vskip 10pt
\begin{defn}\label{def:markedschottkyspace}
The space of equivalence classes of (marked) classical Schottky
groups is called the {\it marked classical Schottky space}; we
denote it by ${\mathcal S}^{\rm mc}_{\rm alg}$.
\end{defn}

Note that we define (marked) classical Schottky groups in terms of
representations rather than the subgroup $\rho(F_n)$ of $\PSLTwoC$
and the space ${\mathcal S}^{\rm mc}_{\rm alg}$ as the space of
such representations, modulo conjugation. In particular, two
representations $\rho_1$ and $\rho_2$ may have the same image in
$\PSLTwoC$ modulo conjugation, but still represent different
points of ${\mathcal S}^{\rm mc}_{\rm alg}$ because of the
marking. To simplify notation, we use $\rho$ to denote elements of
${\mathcal S}^{\rm mc}_{\rm alg}$, instead of $[\rho]$ which is
more cumbersome, there should be no confusion as the traces and
complex lengths are well-defined on the equivalence classes.
\vskip 5pt

Next we give a natural parametrization of the marked classical
Schottky space ${\mathcal S}^{\rm mc}_{\rm alg}$ by the  ideal
fixed points and the ${\rm (trace)}^2$/complex lengths of
$\rho(a_i)$, $i=1, \ldots,n$.

We may normalize $\rho$  by conjugation so that

\centerline{${\rm Fix}^{-}\rho(a_1)=0, {\rm
Fix}^{+}\rho(a_1)=\infty$ and ${\rm Fix}^{-}\rho(a_2)=1$.} \nnn
Then it is not difficult to see that we can parameterize $\rho$ by
\begin{eqnarray*}&&(\,{\rm
Fix}^{+}\rho(a_2), {\rm Fix}^{-}\rho(a_3), {\rm Fix}^{+}\rho(a_3),
\cdots, {\rm Fix}^{+}\rho(a_n); {\rm tr}^2\,
\rho(a_1),\cdots, {\rm tr}^2\,\rho( a_n))\\
&&\in \Cinfty^{2n-3} \times {\mathbf C}^{n},\end{eqnarray*} or,
alternatively, by
\begin{eqnarray*}&&(\,{\rm
Fix}^{+}\rho(a_2), {\rm Fix}^{-}\rho(a_3), {\rm Fix}^{+}\rho(a_3),
\cdots, {\rm Fix}^{+}\rho(a_n); l\,(
\rho(a_1)),\cdots, l(\,\rho( a_n)))\\
&&\in \Cinfty^{2n-3} \times (\CmodTwoPiIZ)^{n},\end{eqnarray*}
where the complex lengths $l(\rho(a))$ are related to the traces
by the formula $$2\cosh \frac{l(\rho(a))}{2}={\rm tr}\,\rho( a),$$
and can be chosen to have positive real part (since all elements
are strictly loxodromic). Then the lengths are defined up to
multiples of $2 \pi i$, and depend only on $|{\rm tr}\,\rho( a)|$
or ${\rm tr}^2\,\rho( a)$. We shall see later that to define the
half-length we will be choosing a lift of the representation and
using the negative of the trace on the right hand side of the
formula. With this normalized parametrization we have

\begin{lem}\label{lem:connected} {\rm (Maskit \cite{maskit2002cm311})}
The marked classical Schottky space ${\mathcal S}^{\rm mc}_{\rm
alg}$ is a path connected open subset of $\Cinfty^{2n-3} \times
(\CmodTwoPiIZ)^{n}$.
\end{lem}

\begin{pf} Here we would like to sketch the idea of the proof;
see Maskit \cite{maskit2002cm311} for a detailed proof. We use the
conformal ball model of hyperbolic 3-space $\HHH$. The ideal
sphere is then the unit sphere $S_{\infty}$. Consider a marked
classical Schottky group $\rho$ so that $\Gamma=\rho(F_n) \subset
{\rm Isom}^{+}(\HHH)$. Then there exists open disks $D_j,
D_j^{\prime}$, $j=1,\cdots,n$ on the ideal sphere whose boundary
circles are denoted $C_j, C_j^{\prime}$ respectively, such that
$\rho(a_j)(C_j)=C_j^{\prime}$ and $\rho(a_j)(D_j)\cap
D_j^{\prime}=\emptyset$. Now we may find a one parameter family
$\rho_t$ of  $\rho$ in ${\mathcal S}^{\rm mc}_{\rm alg}$, such
that the ideal fixed points of all the $\rho_t(a_j)$ are
unchanged, but the real part of $l(\rho_t(a_j)) \to \infty$ (that
is, the translation length of $\rho_t(a_j)$ approaches $\infty$).
We may assume that under the deformation the circles $C_j,
C_j^{\prime}$ shrink towards the respective ideal fixed points of
$\rho(a_j)$, so that their sizes become arbitrarily small. This
gives us room to continue to deform the marked classical Schottky
group in ${\mathcal S}^{\rm mc}_{\rm alg}$ continuously to one
where the fixed points are in some fixed standard configuration,
say, with all of them lying on a circle. Finally, keeping the
ideal fixed points fixed, we can deform the lengths to some
predetermined quantities, with sufficiently large real part. Since
any marked classical Schottky group can be so deformed, the space
${\mathcal S}^{\rm mc}_{\rm alg}$ is path connected. That
${\mathcal S}^{\rm mc}_{\rm alg}$ is open is easily seen from the
definition and parametrization.
\end{pf}

\vskip 5pt

With the parametrization, for each $g \in F_n$, $\rho(g)\in
\PSLTwoC$ is completely determined by the parameters, where $\rho$
is the normalized representation. Hence,  for each $g \in F_n$,
${\rm tr}\,\rho( g)$ and $l(\rho(g))$ are  analytic functions of
the parameters. Furthermore, if we start with a fuchsian Schottky
group, we may define all the lengths to be real and positive, and
if  we extend the  definition of the complex lengths by analytic
continuation on the space $\MCSS$, then the following proposition
states that ${\mathfrak {R}}(l(\rho(g)))>0$ for all $\rho \in
\MCSS$, $g \in F_n$.

\begin{prop}\label{prop:positivereal}
If $\rho_0 \in \MCSS$ is fuchsian, we may define $l(\rho_0(g))$ so
that $l(\rho_0(g))$ is positive real  for all $g \in F_n$. Then if
$l(\rho(g))$ is defined by analytic continuation along a path in
the space $\MCSS$, we have ${\mathfrak R}(l(\rho(g)))>0$ for all
$\rho \in \MCSS$ in the path and all $g \in F_n$.

\end{prop}

\begin{pf}
If $\rho_0$ is fuchsian, then it is well known that $|{\rm
tr}\,\rho _0( g)|>2$ for all $g \in F_n$ so that $l(\rho_0(g))$ is
positive real. If $\rho \in \MCSS$ and $\rho_t$, $0 \le t \le 1$
is a path from $\rho_0$ to $\rho$, then we claim that ${\mathfrak
R}(l(\rho_t(g)))>0$ for all $g \in F_n$, $0 \le t \le 1$.
Otherwise, we will have that ${\mathfrak R}(\rho_t(g))=0$ for some
$g \in F_n$, $t \in [0,1]$, which is impossible as all elements
are strictly loxodromic in a classical Schottky group.
\end{pf}

\vskip 5pt

Finally let us say a few words about the fundamental domain in
$\HHH$ of a classical Schottky group. Let $\Gamma=\rho(F_n)$ be
the image of a classical Schottky group, with the set of disjoint
discs $D_j, D_j^{\prime}$ defined as before. Suppose the circles
$C_j, C_j^{\prime}$ bound respectively geodesic planes
$E_j,E_j^{\prime}$ in $\HHH$. For $j=1,\cdots,n$ let $H_j$ be the
open half space of $\HHH$ bounded by $D_j$ and $E_j$; and
similarly for $H_j^{\prime}$. Then $\mathcal D:=\HHH -
\bigcup_{j=1}^{n}H_j - \bigcup_{j=1}^{n}H_j^{\prime}$ is a
fundamental domain in $\HHH$ of  $\Gamma$. Note that
$D=\overline{\mathcal D} \cap {\mathbf C}_{\infty}$ is the
fundamental domain in the extended complex plane $\Cinfty$ of
$\Gamma$ as described at the beginning of this section. It is
well-known that the quotient hyperbolic 3-manifold $\mathcal H
=\HHH /\Gamma = {\mathcal D} /\Gamma$ is a handlebody of genus
$n$. Note that $\Cinfty /\Gamma = D/\Gamma$ is a conformal surface
of genus $n$ which is the conformal boundary of the handlebody
$\mathcal H$.

\vskip 15pt
\section{McShane's identity for Schottky groups}\label{s:mcshane schottky}
\vskip 5pt

In this section we state and prove our main theorem, Theorem
\ref{thm:mcshane schottky} below. We start with the following
definition.

\vskip 5pt
\begin{defn}\label{def:fuchsianmarking}
A fuchsian marking in ${\mathcal S}^{\rm mc}_{\rm alg}$ is an
element $\rho_0 \in {\mathcal S}^{\rm mc}_{\rm alg}$ which can be
conjugated so that $\rho_0(a_i) \in \PSLTwoR$, for $i=1, \ldots,
n$ (and hence $\rho_0(F_n) \subseteq \PSLTwoR$), and $C_i,C_i'$
can all be taken to be circles which are orthogonal to $\mathbf
R$.

\end{defn}

For a fuchsian marking $\rho_0$, $\HH/\rho_0(F_n)$ is a complete
hyperbolic surface. Its convex core, $M_0$, is a hyperbolic
surface with geodesic boundary, which we call the hyperbolic
surface corresponding to the fuchsian marking. Let $\Delta_0,
\Delta_1, \ldots, \Delta_m$ be the boundary components of $M_0$.
The image $\rho(F_n)$, and hence $F_n$, can be identified with
$\pi_1(M_0)$, and if we define an equivalence relation $\sim$ on
$F_n$ by $g \sim h$ if $g$ is conjugate to $h$ or $h^{-1}$, then
there is a bijection
$${\mathfrak f}:F_n/\sim ~ \to {\mathcal C}$$ from $F_n/\sim$ to the set
${\mathcal C}$ of free homotopy classes of closed curves on $M_0$.
Note that there is a unique geodesic representative for each
non-trivial element of ${\mathcal C}$.

\begin{defn}\label{def:sets}
For a fixed fuchsian marking $\rho_0$, let $M_0$ be the
corresponding hyperbolic surface. Let $\Delta_0, \Delta_1, \ldots,
\Delta_n$ be the boundary components of $M_0$, and let $[d_i]\in
F_n/\sim$, $i=0,\ldots,m$ be the equivalence class corresponding
to the boundary component $\Delta_i$, that is ${\mathfrak
f}[d_i]=\Delta_i$.

We define ${\mathcal P}$ to be the set of all unordered pairs
$\{[g],[h]\}$ of elements in $F_n/\sim$ such that ${\mathfrak
f}[g]$ and ${\mathfrak f}[h]$ bound together with $\Delta_0$ an
embedded pair of pants in $M_0$ (note that it is possible that
${\mathfrak f}[g]=\Delta_k$, for some $1 \le k \le m$).

For $j=1, \ldots, m$, we define ${\mathcal B}_j$ to be the set of
elements $[g]\in F_n/\sim$ such that ${\mathfrak f}[g]$ bounds
together with $\Delta_0$ and $\Delta_j$ an embedded pair of pants
in $M_0$.
\end{defn}

We will also need to define the half lengths, for which we need
representations into $\SLTwoC$ instead of $\PSLTwoC$.

\begin{defn}\label{def:halflength}
If $\rho \in \MCSS$ and $\tilde \rho$ is a lift of $\rho$ to
$\SLTwoC$, then for an element $g \in F_n$, we define the half
length $l(\tilde \rho(g))/2 \in \CmodTwoPiIZ$ of $\tilde \rho (g)$
by
\begin{eqnarray}\label{def:halflength}
\cosh \frac{l(\tilde \rho(g))}{2}=-\frac{{\rm tr}\,\tilde \rho(
g)}{2},
\end{eqnarray} with ${\mathfrak R}l(\rho(g))>0$.
\end{defn}

Note that the real part of the half length is just half of the
real part of the length, and both are positive, while the above
choice fixes the imaginary part, up to multiples of $2 \pi i$. The
minus sign on the right hand side of (\ref{def:halflength}) is
crucial, see Remark \ref{rmk:maintheorem} (iv).

Our main theorem can be stated as follows.

\begin{thm}\label{thm:mcshane schottky}
Let $\rho \in \MCSS$, and let $\tilde \rho$ be any lift of $\rho$
to $\SLTwoC$. Suppose $\rho_0$ is a fuchsian marking, with
corresponding hyperbolic surface $M_0$, with boundary components
$\Delta_0, \ldots, \Delta_m$. Let ${\mathcal P}$ and ${\mathcal
B}_j$, $j=1, \ldots, m$ be defined as in Definition
\ref{def:sets}, relative to $M_0$. Then
\begin{eqnarray}\label{eqn:mcshane schottky}
& &\sum_{\{[g], [h]\}\in {\mathcal P}} \,G\bigg(\frac{l(\tilde
\rho(d_0))}{2},\frac{l(\tilde \rho(g))}{2},
\frac{l(\tilde \rho(h))}{2}\bigg)\nonumber\\
&+& \sum_{j=1}^{m}~\sum_{[g]\in {\mathcal B}_j}
S\bigg(\frac{l({\tilde \rho}(d_0))}{2},\frac{l({\tilde
\rho}(d_j))}{2},\frac{l(\tilde \rho(g))}{2}\bigg)= \frac{l({\tilde
\rho}(d_0))}{2} \mod \pi i.
\end{eqnarray}
 Moreover, each series on the left-hand side of
{\rm(\ref{eqn:mcshane schottky})} converges absolutely.
\end{thm}

\begin{rmk}\label{rmk:maintheorem}\hfill
\begin{enumerate}
\item[(i)] In the case where $\rho=\rho_0$, the above is just a
reformulation of Theorem \ref{thm:complexified non-cusp cases} for
the case of a hyperbolic surface with geodesic boundary
components, and is true without the modulo condition. In fact, the
lift can be chosen so that the right hand side is real and
positive.

\item[(ii)] The identity (\ref{eqn:mcshane schottky}) is true only
modulo $\pi i$ because we have fixed the choice of the
$\tanh^{-1}$ function in the definition of the functions
$G(x,y,z)$ and $S(x,y,z)$ (definition \ref{def:GandS}), this may
differ from the values obtained by analytic continuation by some
multiple of $\pi i$. Indeed, as we will see in Corollary
\ref{cor:mcshane schottky m=0} and in the example in \S
\ref{s:example schottky}, there is a difference of $2 \pi i$ in
that example, where $m=0$, that is, $M_0$ has a single boundary
component.

\item[(iii)] The result is independent of the lift chosen. This is
because if $\tilde \rho$ and $\bar \rho$ are two different lifts
of $\rho$, then for each of the summands on the first series,
either ${\rm tr}\,\tilde \rho( g), {\rm tr}\,\tilde \rho( h)$ and
${\rm tr}\,\tilde \rho( d_0)$ are all equal to ${\rm tr}\,\bar
\rho( g), {\rm tr}\,\bar \rho( h)$ and ${\rm tr}\,\bar \rho( d_0)$
or exactly two of them differ by their signs (and similarly for
the summands in the second series). In the latter case, two of the
half lengths differ by $\pi i$, but it can be easily checked that
both $G(x,y,z)$ and $S(x,y,z)$ remain the same if $\pi i$ is added
to two of the arguments.

\item[(iv)] The choice of the half length functions given above is
not arbitrary but arises from the computation of $G(x,y,z)$ and
$S(x,y,z)$ as ``gap'' functions (this is based on the convention
adopted by Fenchel in \cite{fenchel1989book}, see
\cite{tan-wong-zhang2004cone-surfaces} for details; see also
\cite{goldman2003gt} where Goldman uses a similar convention).
Roughly speaking, the relative positions of the axes for $\tilde
\rho( g)$, $\tilde \rho( h)$ and $\tilde \rho( d_0)$ are
completely determined by their traces. These axes form the
non-adjacent sides of a right angled hexagon in $\HHH$ and the
half lengths basically arise as the lengths of these sides of the
hexagon.

\item[(v)] It can be shown that $G(x,y,z)$ and $S(x,y,z)$ can also
be expressed as
\begin{eqnarray}\label{eqn:G in log}
G(x,y,z)=\log\frac{\exp(x)+\exp(y+z)}{\exp(-x)+\exp(y+z)},
\end{eqnarray}
\begin{eqnarray}\label{eqn:S in log}
S(x,y,z)=\frac{1}{2}\log\frac{\cosh(z)+\cosh(x+y)}{\cosh(z)+\cosh(x-y)},
\end{eqnarray}
as used by Mirzakhani in \cite{mirzakhani2004preprint} (with
different notation), where the function $\log$ is in the principal
branch, that is, the imaginary parts of its images are in $(-\pi,
\pi]$. It would be interesting to see if her results on the
Weil-Petersson volumes can be generalized to the classical
Schottky space.
\end{enumerate}

\end{rmk}

\begin{pf} Let $\rho_t$, $0\le t \le 1$ be a
deformation from the fuchsian marking to an arbitrary marked
classical Schottky group $\rho$, where $\rho_0$ is the fuchsian
marking and $\rho_1=\rho$; this is possible by Lemma
\ref{lem:connected}. Let $\tilde \rho_t$ be a continuous lift of
$\rho_t$. We shall then prove that the series on the left-hand
side of (\ref{eqn:mcshane schottky}) converge uniformly on compact
subsets of the marked classical Schottky space ${\mathcal S}_{\rm
alg}^{\rm mc}$. Then each side of (\ref{eqn:mcshane schottky}) is
a holomorphic function modulo $\pi i$ on the space ${\mathcal
S}_{\rm alg}^{\rm mc}$. By Theorem \ref{thm:complexified non-cusp
cases} the identity (\ref{eqn:mcshane schottky}) holds on the
totally real subspace in a neighborhood of $\rho_0$ in ${\mathcal
S}_{\rm alg}^{\rm mc}$. Note that for this to be true, the correct
choice of the half length as given in (\ref{def:halflength}) must
be used, see \cite{tan-wong-zhang2004cone-surfaces} for details.
 Hence the identity also holds modulo $\pi i$ for each
$t \in [0,1]$ by analytic continuation. This proves the identity
for a particular lift; that it holds for all lifts now follow from
Remark \ref{rmk:maintheorem}(iii).

Given a compact subset ${\mathcal K}$ of ${\mathcal S}_{\rm
alg}^{\rm mc}$, we have a constant $\kappa >0$ such that for each
 $\rho \in {\mathcal K}$, $\rho(F_n)$ has a
fundamental domain $\mathcal D$ in $\HHH$ as described at the end
of \S \ref{s:intro schottky}, and the minimum hyperbolic distance
between any pair of its bounding geodesic planes
$E_1,E_1^{\prime}, \cdots, E_n, E_n^{\prime}$ is $\ge \kappa$.
Then we have the following length estimate lemma for $\rho \in
{\mathcal K}$.

\begin{lem}\label{lem:wordlengthbound} If $g \in F_n$ is a cyclically
reduced word in the set of generators $a_1,
a_1^{-1}, \cdots, a_n, a_n^{-1}$ with word length $\| g\| $, then
the closed geodesic $\gamma$ which $\rho(g)$ represents in the
quotient hyperbolic $3$-manifold $\mathcal D=\HHH/\rho(F_n)$ has
hyperbolic length $\ge \kappa \| g\| $.
\end{lem}

\begin{pf} Choose in $\HHH$ an arbitrary lift, $\tilde{\gamma}$,
of the closed geodesic $\gamma$. Note that $\HHH$ is tiled by the
images of a fundamental domain $\mathcal D$ under the action of
elements of $\rho(F_n)$, that is, $\HHH=\bigcup_{g'\in
F_n}\rho(g')(\mathcal D)$. It can be shown that the line
$\tilde{\gamma}$ in $\HHH$ passes through $\| g\| $ successive
images of $\mathcal D$ ``periodically'' dictated by the word $g$.
Thus the hyperbolic length of $\gamma$, which equals the length of
the part of $\tilde{\gamma}$ lying in the union of these $\| g\| $
successive images of $\mathcal D$, is at least $\kappa \,\| g\| $.
\end{pf}

Now we prove the uniform convergence of the first series in
(\ref{eqn:mcshane schottky}) for $t \in [0,1]$. By the above lemma
there is a constant $\kappa >0$ such that for every $t \in [0,1]$
and every $g\in F_n$, we have $L(\tilde \rho_t(g)) \ge \kappa \|
g\| $, where $L(\tilde \rho_t(g))$ is the hyperbolic length of the
closed geodesic that $\tilde \rho_t(g)$ represents in the quotient
hyperbolic 3-manifold $\HHH/\tilde \rho_t(F_n)$, and where $\| g\|
$ is the cyclically reduced word length of $g$ in the letters
$a_1^{\pm 1},\cdots, a_n^{\pm 1}$.

Note that the image of the fuchsian marking $\rho_0(F_n) \subset
\PSLtwoR$ has a fundamental domain $\mathcal D(0)$ in $\HHH$ whose
intersection with $\HH \subset \HHH$ is a fundamental domain  of
$\HH/\rho_0(F_n)$ in $\HH$.  Let $\mathcal P$ be the set of
unordered pairs of elements $\{[g],[h]\}$ of $F_n/\sim$ as defined
in Definition \ref{def:sets}. Then the pairs $\{[g],[h]\} \in
{\mathcal P}$ can be ordered by using the sum of their (cyclically
reduced) word lengths $\| g\|+\|h\|$.  We have the following
version of the Birman-Series result:

\begin{lem}\label{lem:Birman-Series}{\rm (c.f. Lemma 2.2
\cite{birman-series1985t})}\, There exists a polynomial  $P(n)$
such that the number of pairs $\{[g], [h]\}$ in $\mathcal P$ with
$\| g\| +\| h\| =n$ is no greater than $P(n)$.
\end{lem}

\begin{pf}
The proof follows the same line as that used in
\cite{birman-series1985t}. We start with the fundamental polygon
$\mathcal D(0) \cap \HH$, then for each pair  $\{[g], [h]\}$ in
$\mathcal P$ with $\| g\| +\| h\| =n$ we can associate a simple
diagram on $\mathcal D(0) \cap \HH$ consisting of $n$ disjoint
arcs (note that the original Birman--Series' argument in
\cite{birman-series1985t} is for just one simple geodesic, but it
works as well here for the pair of disjoint simple closed
geodesics on the surface $M_0$ associated with $\{[g], [h]\}$).
Conversely, the diagram determines the pair $\{[g], [h]\}$). The
number of such simple diagrams is bounded by a polynomial $P(n)$.
\end{pf}

Now, for $t \in [0,1]$ and $[g]\in F_n/\sim$, the real part $\Re\,
l(\tilde \rho_t(g))$ of the complex translation length $l(\tilde
\rho_t(g))$ is equal to the hyperbolic length $L(\gamma_t(g))$ of
the closed geodesic $\gamma_t(g)$, which represents $[g]$ in the
quotient hyperbolic manifold $\HHH/\tilde \rho_t(F_n)$. Hence for
a pair $\{[g], [h]\}$ in $\mathcal P$, where  $g,h$ are cyclically
reduced, we have
\begin{eqnarray}
\bigg|G\bigg(\frac{l(\tilde \rho_t(d_0))}{2},\frac{l(\tilde
\rho_t(g))}{2},\frac{l(\tilde \rho_t(h))}{2}\bigg)
\bigg|\nonumber&\le& {\rm const} \cdot \exp \bigg(-
\frac{L(\gamma_t(g))+L(\gamma_t(h))}{2} \bigg)\nonumber\\ &\le&
{\rm const} \cdot \exp \bigg(- \frac{\kappa(\| g\| +\| h\| )}{2}
\bigg),
\end{eqnarray}
where the first inequality follows from (\ref{eqn:G in log}) and
the fact that $|\log(1+u)|\le 2|u|$ for log in the principal
branch and for $u \in {\mathbf C}$ with $|u| \le \frac12$, and the
the last inequality follows from Lemma \ref{lem:wordlengthbound}.
It now follows Lemma \ref{lem:Birman-Series} that the first series
in (\ref{eqn:mcshane schottky}) converges absolutely and uniformly
for $t \in [0,1]$.

The absolute and uniform convergence for the other series in
(\ref{eqn:mcshane schottky}) can be similarly proved.
\end{pf}

\begin{cor}\label{cor:mcshane schottky m=0}
If $m=0$ in Theorem {\rm \ref{thm:mcshane schottky}} then we have
\begin{eqnarray}\label{eqn:mcshane schottky m=0}
& &\sum_{\{[g], [h]\}\in {\mathcal P}} \,G\bigg(\frac{l(\tilde
\rho(d_0))}{2},\frac{l(\tilde \rho(g))}{2}, \frac{l(\tilde
\rho(h))}{2}\bigg) = \frac{l(\tilde \rho(d_0))}{2} \mod 2 \pi i,
\end{eqnarray}
{\rm(}note that here the identity holds modulo $2 \pi i$ instead
of $\pi i${\rm)} and the series converges absolutely. \square
\end{cor}

\vskip 3pt

Note that if $m=0$, then $d_0$ is actually a commutator so that
${\rm tr}(\tilde \rho(d_0))$ is independent of the lift $\tilde
\rho$, hence,  the right hand side of (\ref{eqn:mcshane schottky
m=0}) is independent of the lift. In fact, it can be shown, with a
little bit of extra work, that for the fuchsian marking $\rho_0$,
${\rm tr}(\tilde \rho_0(d_0))$ is always strictly negative in this
case (see for example \cite{goldman2003gt}). Furthermore, for the
more general case where $m \neq 0$, we can always choose a lift
$\tilde \rho_0$ such that ${\rm tr}(\tilde \rho_0(d_0))$ is
strictly negative.

\vskip 15pt
\section{The Weierstrass Indentities}\label{s:Weierstrass}
\vskip 5pt

In this section, we consider rank 2 classical Schottky groups
$\rho:F_2 \to \PSLTwoC$. For ease of notation, we denote the
marked generators of $F_2$ by $a$ and $b$. Let $\rho_0$ be a
fuchsian marking such that the corresponding surface $M_0$ is the
one-holed torus with geodesic boundary $\Delta_0$.  Let ${\mathcal
S}$ be the set of non-peripheral simple closed geodesics on $M_0$
and let $w_1,w_2$ and $w_3$ be the Weierstrass points on $M_0$.
Then each element of ${\mathcal S}$ passes through exactly two of
the Weierstrass points, and we define the Weierstrass classes to
be the subsets ${\mathcal A}_i$, $i=1,2,3$ of ${\mathcal S}$
consisting of those geodesics which miss $w_i$. Then ${\mathcal
S}=\bigsqcup_{i=1}^3 {\mathcal A}_i$. Let $$\bar {\mathcal
S}:=\{[g] \in F_2/\sim \,| \,{\mathfrak f}[g] \in {\mathcal S}\},
\quad \bar {\mathcal A}_i:=\{[g]\in F_2/\sim \,| \,{\mathfrak
f}[g] \in {\mathcal A_i}\}$$ be the corresponding sets in
$F_2/\sim$. Then $[g] \in \bar {\mathcal S}$ if and only if any
cyclically reduced representative $g$ forms with another element
$h$ a generating set for $F_2$. The set $\bar {\mathcal S}$  can
be identified with ${\mathbf Q}\cup \infty$ by considering the
slopes of the corresponding simple closed curves on the torus
(without the hole), and the subsets $\bar {\mathcal A}_i$,
$i=1,2,3$ with the subsets of ${\mathbf Q}\cup \infty$ with both
numerator and denominator odd, numerator odd and denominator even,
and numerator even and denominator odd, respectively; see
\cite{tan-wong-zhang2004gMm} for details.

 We have the following extension of the Weierstrass identities
proven by McShane in \cite{mcshane2004blms} (see also
\cite{tan-wong-zhang2004cone-surfaces}).

\begin{thm}\label{thm:Weierstrass}
For any rank 2 classical Schottky group $\rho \in \MCSS$, if
$\bar{\mathcal A}$ is a Weierstrass class, and $\tilde \rho$ is a
lift of $\rho$, then
\begin{eqnarray}\label{eqn:weierstrass}
\sum_{[g] \in \bar{\mathcal A}} \tan^{-1} \left ( \frac{\cosh
\frac{l(\tilde \rho(d_0))}{4}}{\sinh \frac{l(\tilde \rho(g))}{2}}
\right )= \frac{\pi}{2} \mod \pi,
\end{eqnarray}
where $[d_0]=[b^{-1}a^{-1}ba]\in F_2/\sim$ corresponds to the
boundary $\Delta_0$ of $M_0$.

\end{thm}
Note that there are two choices for the quarter-length
$\frac{l(\tilde \rho(d_0))}{4}$, we can choose either, but the
choice should be the same for all summands.  The half length
$\frac{l(\tilde \rho(g))}{2}$ depends on the choice of the lift
$\rho$.

We skip the proof as it is essentially the same as the proof for
Theorem \ref{thm:mcshane schottky}. Here is a geometric
interpretation for the above result. First note that the case of
the two generator group is very special since for a lift $\tilde
\rho$, each of $\tilde \rho(a_1)$ and $\tilde \rho(a_2)$ can be
factored as a product of two involutions (half turns) as follows:
$$\tilde \rho(a_1)=-H_3H_2, \qquad \tilde \rho(a_2)=-H_1H_3$$
where $H_i \in \SLTwoC$ with $H_i^2=-I$ and the axis $l_i$ for
$H_i$ are such that $l_2$ and $l_3$ are perpendicular to the axis
for $\rho(a_1)$ and $l_3$ and $l_1$ are perpendicular to the axis
for $\rho(a_2)$ (the minus sign is the convention adopted in
Fenchel \cite{fenchel1989book}). Furthermore, in the case of the
fuchsian marking $\rho_0$, $l_i$ are lines in $\HHH$ perpendicular
to $\HH$ and passing through lifts of the Weierstrass points $w_i$
respectively. Calling $l_i$ the Weierstrass axes, a deformation
$\rho_t$ will correspond to a deformation of the relative
positions of the Weierstrass axes. The summands of the left hand
side of (\ref{eqn:weierstrass}) are then gaps arising from the
fixed points of $\rho(g)$ (for certain representatives $g$ of
$[g]$) measured against the Weierstrass axis $l$ corresponding to
the Weierstrass class $\bar{\mathcal A}$.

There are also generalizations of the variations and refinements
of McShane's identity given by Bowditch in \cite{bowditch1997t},
Sakuma in \cite{sakuma1999sk}, Akiyoshi, Miyachi and Sakuma in
\cite{akiyoshi-miyachi-sakuma2002preprint},
\cite{akiyoshi-miyachi-sakuma2004preprint}. In those cases, the
variations of McShane's identity were defined relative to a cusp,
these can be generalized to the case of identities relative to a
boundary geodesic. We can then study deformations (say in the
space of discrete, faithful representations) where the trace of
the boundary component remains real and with absolute value $>2$.
Most of the identities can then be re-interpreted and generalized
to this context.

 \vskip 15pt
\section{An example: The three-holed sphere}\label{s:example schottky}
 \vskip 5pt

As in the previous section, we consider rank 2 classical Schottky
groups $\rho \in \MCSS$, with the same notation. In particular,
$F_2=\langle a,b\rangle$; $\rho_0$ is a fuchsian marking of
$\MCSS$ where the corresponding hyperbolic surface $M_0$ is a
one-holed torus with geodesic boundary $\Delta_0$; and $\bar
{\mathcal S}\subset F_2/\sim$ consists of the equivalence classes
corresponding to the non-peripheral simple closed geodesics on
$M_0$. Let $\rho_1 \in \MCSS$ correspond to another fuchsian
marking, with corresponding hyperbolic surface $M_1$, a pair of
pants with geodesic boundary, and let $\hat {\mathcal S}$ be the
set of closed geodesics on $M_1$ corresponding to the elements of
$\bar {\mathcal S}$. Note that the three boundary geodesics of
$M_1$ are elements of $\hat{\mathcal S}$ corresponding to $[a],
[b]$ and $[ab]$ (with appropriate orientation). Apart from these,
all other elements of $\hat{\mathcal S}$ are non-simple geodesics
on $M_1$.

\begin{figure}
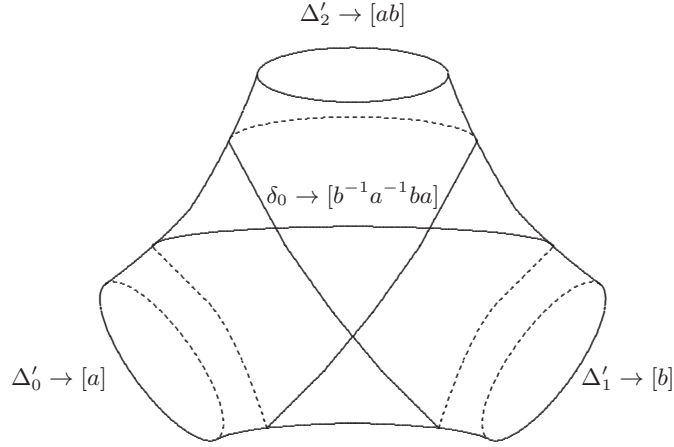

\begin{center}
\mbox{\beginpicture \setcoordinatesystem units <0.05in,0.05in>
\setplotarea x from -40 to 40 , y from 0 to 45

\ellipticalarc axes ratio 3.5:1 360 degrees from 10 40  center at
0 40

\startrotation by .809 .588 about -20 10 \ellipticalarc axes ratio
1:3 180 degrees from -20 20 center at -20 10 \stoprotation

\startrotation by .809 -.588 about 20 10 \ellipticalarc axes ratio
1:3 -180 degrees from 20 20 center at 20 10 \stoprotation

\startrotation by .809 .588 about -20 10 \setdashes<1.50pt>
\ellipticalarc axes ratio 1:3 -180 degrees from -20 20 center at
-20 10 \stoprotation

\startrotation by .809 -.588 about 20 10 \ellipticalarc axes ratio
1:3 180 degrees from 20 20 center at 20 10 \stoprotation \setsolid

\plot -9 3 -3 9 2 15 7 22 13 33 /

\plot 9 3 3 9 -2 15 -7 22 -13 33 /

\ellipticalarc axes ratio 10:1 180 degrees from 21 22  center at 0
22

\setquadratic \plot -25.88 18.09 -21 22 -17 26 -13 33 -10 40 /

\setquadratic \plot 25.88 18.09 21 22 17 26 13 33 10 40 /

\setquadratic \plot -14.12 1.89 -9 3 0 3.5 9 3 14.12 1.89 /

\setdashes<1.50pt>

\ellipticalarc axes ratio 5:1 180 degrees from 13 33  center at 0
33

\plot -9 3 -12 11 -15 16 -18  19  -21 22 /

\plot 9 3 12 11 15 16 18  19  21 22 /\setsolid

\put {\mbox{\small $\Delta_0' \to [a]$}} [cb] <-1mm,-4mm> at -30
10

\put {\mbox{\small $\Delta_1' \to [b]$}} [cb] <1mm,-4mm> at 28 10

\put {\mbox{\small $\Delta_2' \to [ab]$}} [cb] <0mm,0mm> at 0 45

\put {\mbox{\small $\delta_0 \to [b^{-1}a^{-1}ba]$}} [cb]
<0mm,0mm> at 0 26
\endpicture}
\end{center}
\caption{A commutator curve on $M_1$}\label{fig:comm}
\end{figure}

Let $\delta_0$ be the geodesic on $M_1$ corresponding to the
commutator  $[d_0]=[b^{-1}a^{-1}ba]$. Note that $\delta_0$ is not
a simple geodesic on $M_1$, in fact it has triple
self-intersection; see Figure \ref{fig:comm}.

Let the three geodesic boundary components of $M_1$ be denoted by
$\Delta_0', \Delta_1', \Delta_2'$ with hyperbolic lengths
$l_0=l(\Delta_0') > 0, l_1=l(\Delta_1') > 0, l_2=l(\Delta_2')
> 0$ respectively. Then for this hyperbolic surface, we have the
following trivial identity
\begin{eqnarray}\label{eqn:3p trivial}
G\bigg(\frac{l_0}{2},\frac{l_1}{2},\frac{l_2}{2}
\bigg)+S\bigg(\frac{l_0}{2},\frac{l_1}{2},\frac{l_2}{2}\bigg)+
S\bigg(\frac{l_0}{2},\frac{l_2}{2},\frac{l_1}{2}\bigg)=\frac{l_0}{2}.
\end{eqnarray}

\vskip 5pt

There is, however, a non-trivial identity on $M_1$ derived from
the fuchsian marking $\rho_0$. Recall that the trace  ${\rm
tr}\,\rho[d_0]$ is well-defined and independent of the lift of
$\rho$, since $d_0$ is a commutator. In fact, for $\rho_0$, we
have ${\rm tr}\,\rho_0[d_0]<-2$, and for $\rho_1$, we have ${\rm
tr}\,\rho_1[d_0]>18$ (see Goldman \cite{goldman2003gt} for details
where he studied geometric structures arising from 2 generator
subgroups of $\PSLTwoC$ with real character varieties in detail).
In particular, the half length $l(\delta_0)/2$ of $\delta_0$ is
well-defined up to multiples of $2 \pi i$, and
$l(\delta_0)/2=|\delta_0|/2+\pi i$, where $|\delta_0|$ is the
usual hyperbolic length of $\delta_0$ on $M_1$, since ${\rm
tr}\,\rho_1[d_0]>2$ (recall the definition of the half length from
(\ref{def:halflength})). From Corollary \ref{cor:mcshane schottky
m=0}, we have
\begin{eqnarray}\label{eqn:example}
& &\sum_{\alpha \in \hat{\mathcal S}}
\,G\bigg(\frac{l(\delta_0)}{2},\frac{l(\alpha)}{2},
\frac{l(\alpha)}{2}\bigg) = \frac{l(\delta_0)}{2} \mod 2 \pi i,
\end{eqnarray}
or equivalently,
\begin{eqnarray}\label{eqn:schottky hole torus}
\sum_{\alpha \in \hat{\mathcal S}} 2 \tanh^{-1}\bigg( \frac {\sinh
(l(\delta_0)/2)} {\cosh (l(\delta_0)/2)+ \exp l(\alpha)} \bigg)=
\frac{l(\delta_0)}{2} \mod 2 \pi i.
\end{eqnarray}

\vskip 5pt

We see from the second expression that it does not matter which
choice of the two half lengths we use for $\alpha$. It turns out
that the three summands in the left hand side of
(\ref{eqn:example}) or (\ref{eqn:schottky hole torus})
corresponding to the three boundary components of $M_1$ have
positive real part and imaginary part equal to $\pi i$ whereas all
the other summands are real and $<0$. As pointed out earlier, the
right hand side is a complex number with positive real part and
imaginary part equal to $\pi i$.

\begin{figure}
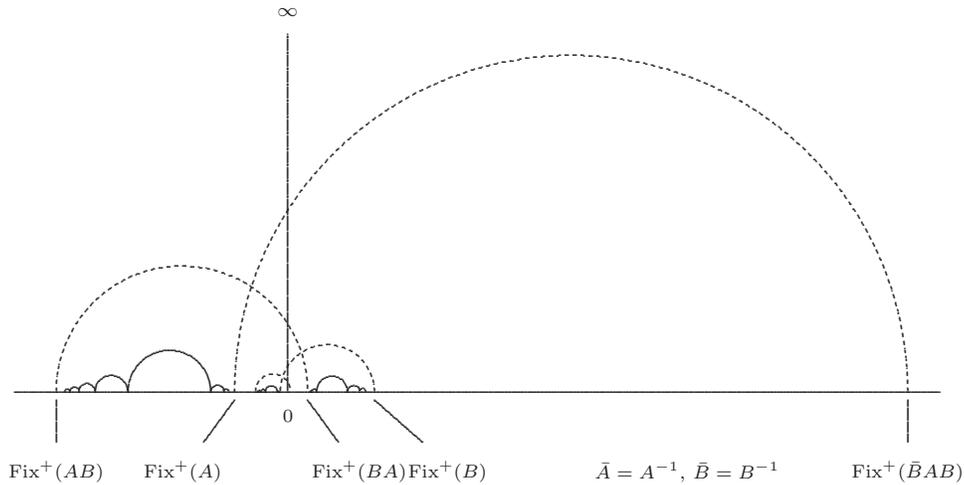


\begin{center}
\mbox{
\beginpicture
\setcoordinatesystem units <0.22in,0.22in>

\setplotarea x from -6.5 to 15.5, y from 0 to 9

\plot -6.5 0 15.5 0 / \plot 0 0 0 8.5 /

\circulararc 180 degrees from -1.520 0  center at -1.674 0

\circulararc 180 degrees from -1.835 0  center at -2.816 0

\circulararc 180 degrees from -3.814 0  center at -4.199 0

\circulararc -180 degrees from 0.568 0  center at 0.625 0

\circulararc -180 degrees from 0.686 0  center at 1.052 0

\circulararc -180 degrees from 1.425 0  center at 1.569 0

\circulararc -180 degrees from -0.188 0  center at -0.192 0

\circulararc -180 degrees from -.196 0  center at -.204 0

\circulararc 180 degrees from -.212 0  center at -.234 0

\circulararc 180 degrees from -.256 0  center at -.393 0

\circulararc 180 degrees from -.532 0  center at -.586 0

\circulararc 180 degrees from -.640 0  center at -.667 0

\circulararc 180 degrees from -.694 0  center at -.709 0

\circulararc 180 degrees from -1.402 0  center at -1.461 0

\circulararc 180 degrees from -1.827 0  center at -1.831 0

\circulararc 180 degrees from -4.587 0  center at -4.778 0

\circulararc -180 degrees from .524 0  center at .546 0

\circulararc -180 degrees from 1.714 0  center at 1.785 0

\circulararc 180 degrees from -4.970 0  center at -5.076 0

\circulararc 180 degrees from -5.182 0  center at -5.245 0



\put {\mbox{\scriptsize $0$}} [cb] <0mm,-4mm> at 0 0

\put {\mbox{\scriptsize $\infty$}} [cb] <0mm,2.5mm> at 0 8.5

\put {\mbox{\scriptsize ${\rm Fix}^{+}({\bar B}AB)$}} [cb]
<0mm,-12mm> at 14.745 0

\plot 14.745 -1.20 14.745 -0.20 /

\put {\mbox{\scriptsize ${\rm Fix}^{+}(A)$}} [cb] <0mm,-12mm> at
-2.5 0

\plot -2.0 -1.20 -1.265 -0.20 /

\put {\mbox{\scriptsize ${\rm Fix}^{+}(B)$}} [cb] <0mm,-12mm> at
3.8 0

\plot 3.2 -1.20 2.059 -0.20 /

\put {\mbox{\scriptsize ${\rm Fix}^{+}(AB)$}} [cb] <0mm,-12mm> at
-5.509 0

\plot -5.509 -1.20 -5.509 -0.20 /

\put {\mbox{\scriptsize ${\rm Fix}^{+}(BA)$}} [cb] <0mm,-12mm> at
1.7 0

\plot 1.2 -1.20 0.47 -0.20 /

\put {\mbox{\scriptsize ${\bar A}=A^{-1},\,{\bar B}=B^{-1}$}} [cb]
<0mm,-12mm> at 9.5 0

\setdashes<1.50pt>

\circulararc 180 degrees from 14.745 0  center at 6.740 0

\circulararc -180 degrees from -5.509 0  center at -2.518 0

\circulararc 180 degrees from 2.058 0  center at 0.941 0

\circulararc -180 degrees from -.769 0  center at -.351 0

\endpicture
}\end{center}

\caption{Gaps for $M_1$}\label{fig:schottky gaps}
\end{figure}

If we define (as in \cite{tan-wong-zhang2004gMm})
$\nu:=l(\delta_0)/2=\cosh^{-1}(-{\rm tr}\, (\rho_1(d_0))/2$, and
the functions $h(x)$ and ${\mathfrak h}(x)$ by
\begin{eqnarray}\label{frak h(x)=}
h(x)=\frac{1}{2}\bigg(1-\sqrt{1-\frac{4}{x^2}}\bigg) \quad {\rm
and} \quad {\mathfrak h}(x)=\log
\bigg(\frac{1+(e^{\nu}-1)h(x)}{1+(e^{-\nu}-1)h(x)}\bigg),
\end{eqnarray}
where the square root is always chosen to have positive real part
and we use the principal branch (with imaginary part in
$(-\pi,\pi]$) for the $\log$ function, then we can also express
(\ref{eqn:schottky hole torus}) as
\begin{eqnarray}\label{eqn:3p frak}
\sum_{[g] \in \bar{\mathcal S}} {\mathfrak h}\big({\rm
tr}\,\rho_1(g)\big)= \nu \mod 2 \pi i.
\end{eqnarray}
\begin{rmk}{\rm The  identity (\ref{eqn:3p frak}) was derived in \cite{tan-wong-zhang2004gMm}
with a different proof and under much more general conditions, and
expressed in terms of  the $\mu$-Markoff map corresponding to the
$\mu$-Markoff triple $(x,y,z)=({\rm tr}\,a , {\rm tr}\,b, {\rm
tr}\,ab)$.} Necessary and sufficient conditions for the identity
to hold were also given in \cite{tan-wong-zhang2004nsc}. For
example, the identity still holds if some or all of the boundary
components of $M_1$ degenerate to cusps.
\end{rmk}
We give a brief description of the geometric interpretation of the
above; see \cite{zhang2004thesis} or \cite{tan-wong-zhang2004gMm}
for more details. For each $[g] \in \bar{\mathcal S}$, we may
choose representatives $g_1$ and $g_2$ such that $g_1g_2=d_0$.
Then each of the summands in the various versions of the identity
above can be interpreted as gaps from the attractive fixed point
of $\rho_1(g_1)$ to the repelling fixed point of $\rho_1(g_2)$
measured along the (directed) axis of $\rho_1(d_0)$. In the cases
where $[g] \in \bar {\mathcal S}$ corresponds to the three
boundary components of $M_1$ (that is, $[a]$, $[b]$ and $[ab]$),
the above fixed points lie on the two different intervals of
${\mathbf R} \cup \infty$ separated by the fixed points of
$\rho_1(d_0)$, which is why the summand has imaginary part $\pi i$
in these three cases. For the other summands, the fixed points lie
on the same interval and the summands are real. See Figure
\ref{fig:schottky gaps} where  we have $\rho_1(a)=A, \rho_1(b)=B$,
where $A, B \in \PSLTwoR$, and we use the notation $\bar
A:=A^{-1}$, $\bar B:=B^{-1}$. The picture is normalized so that
the fixed points of the commutator $\bar B \bar A BA$ are 0 and
$\infty$. In the figure, the gap arising from $[a]\in
\bar{\mathcal S}$ corresponding to one of the boundary components
of $M_1$ is measured by dropping perpendiculars from ${\rm
Fix}^{+}(A)$ and ${\rm Fix}^{+}(\bar B A B)$ to the axis $[0,
\infty]$ of $\bar B \bar A BA$; hence this gap has positive real
part and imaginary part $\pi i$. Note that $(\bar B A
B)^{-1}A=\rho_1(d_0)$. Similarly the other two middle-sized dotted
semi-circles in the figure show the gaps arising from $[b]$ and
$[ab]$ corresponding to the other two boundary components. Gaps
for the other elements of $\bar{\mathcal S}$ can be similarly
obtained by a recursive process using the Farey construction of
the rationals; all these gaps will be real and negative (they are
represented by the solid semi-circles in the figure).

\vskip 5pt

 Finally, we note that more generally, if the rank $n \ge 3$, it is
possible to have two fuchsian markings $\rho_0$ and $\rho_1$ such
that the corresponding hyperbolic surfaces $M_0$ and $M_1$ are
homeomorphic but the markings are different. Then $\rho_0$ will
induce identities for $M_1$ which are different from those
obtained from $\rho_1$; for example the gaps may be measured
against a geodesic on $M_1$ which is not necessarily a boundary
geodesic, or even simple.

\vskip 10pt

{}

\end{document}